\documentclass{article}
\usepackage[utf8]{inputenc}

\title{Stratified formulas are equivalent to Acyclic formulas: A review}
\author{Zuhair Al-Johar}
\date{ 9 February 2018}

\begin{document}

\maketitle

\section{Introduction}
This document include my comments on Nathan Bowler's proof that every stratified formula is equivalent to an acyclic formula. The whole proof of Nathan's is quoted exactly as it is presented in his messages to me. The comments are edits of clarifying e-mail messages that I've sent to Thomas Forster, as well as e-mail correspondence with Nathan Bowler. Minor technical issues are addressed also.

\section {Nathan's proof} 

``Any stratified formula can be rephrased as an acyclic formula (assuming we have enough comprehension to build finite sets $\{x_1, \cdots , x_k\}$ given the elements $x_1, \cdots , x_k$)."

``Basic things we need: let $\iota$ be the operator $x \mapsto \{x\}$
There is an acyclic
expressing formula expressing $y = \iota (x)$ , namely $(\exists w)w \in y \wedge (\forall z)z \in y \to z = x$. Thus
for any natural number $n$ (of the metalanguage) there is an acyclic formula
expressing $y = \iota^n(x)$. There is also an acyclic formula expressing $x = \emptyset$, namely
$(\forall y)y \not \in x$ . Thus for any natural number $n$ there is an acyclic formula expressing
$x = \iota^n(\emptyset)$.

In our paper, under the heading `An aside: Implementation of the Wiener
pair', it is shown that each of `$p$ is a Wiener pair', `$x = \pi_1(p)$' and `$x = \pi_2(p)$'
is expressible with an acyclic formula. I'll express the first as $P(p)$.
We'll be working with functions whose domain is $\{\iota^i(\emptyset)|1 \leq i \leq n\}$ for some
natural number $n$ (of the metalanguage). First, let's check that there is an
acyclic formula (which I'll denote $F_n(f)$), expressing that $f$ is such a function. We can use the formula\bigskip

 $(\forall p) p \in f \to [P(p) \wedge (\pi_1(p)=\iota (\emptyset) \lor \pi_1(p)=\iota^2(\emptyset) \lor \cdots \lor \pi_1(p)= \iota^n(\emptyset))] \wedge (\exists p_1) (p_1 \in f \wedge \pi_1(p_1) = \iota (\emptyset)) \wedge \cdots \wedge(\exists p_n) (p_n \in f \wedge \pi_1(p_n) = \iota^n(\emptyset)) \wedge (\exists y_1)(\forall q_1)[q_1 \in f \wedge \pi_1(q_1)=\iota(\emptyset)] \to \pi_2(q_1)=y_1 \wedge \cdots \wedge (\exists y_n)(\forall q_n)[q_n \in f \wedge \pi_1(q_n)=\iota^n(\emptyset)] \to \pi_2(q_n)=y_n$\bigskip

We can express $f(\iota^i(\emptyset))=x$ with the acyclic formula $(\exists p \in f) (\pi_1(p) = \iota^i(\emptyset) \wedge \pi_2(p)=x)$. We can express $f(\iota^i(\emptyset))=f(\iota^j(\emptyset))$ with the formula $(\exists y)(\forall p) [p \in f \wedge (\pi_1(p) = \iota ^i (\emptyset) \lor \pi_1(p)= \iota^j(\emptyset))] \to \pi_2(p) = y$.  Finally, we need an acyclic
formula which, assuming that $f(\iota^i(\emptyset))$ is of the form $\iota ^{d+1} (x)$, expresses that
there is some y with $f(\iota^j(\emptyset)) = \iota^d(y)$ and $x \in y$. For this we can use the
formula
$$(\exists z)(\forall p)(\forall w)[p \in f \wedge \pi_2(p) = \iota^d(w) \wedge  (\pi_1(p) = \iota^i(\emptyset) \lor \pi_1(p) = \iota^j(\emptyset))] \to z \in w$$
and we denote this by $f(\iota^i(\emptyset)) \in_d f(\iota^j(\emptyset))$.

Now suppose we have a stratified formula $\phi$ for which we want to find an
acyclic equivalent. Because of the existence of prenex normal forms, we may
assume without loss of generality that $\phi$ is quantifier free. Let the variables
appearing in $\phi$ be $x_1,..,x_n$, where the stratification level of $x_i$ is $d_i$. Without loss
of generality all the $d_i$ are negative. Let $\hat{\phi }$  be the acyclic formula in the single free
variable $f$ obtained from $\phi$ by replacing atomic subformulae of the form $x_i = x_j$
by acyclic formulae expressing $f(\iota^i(\emptyset)) = f(\iota^j(\emptyset))$ and atomic subformulae of
the form $x_i \in x_j$ by acyclic formulae expressing $f(\iota^ i(\emptyset)) \in_{-d_j} f(\iota^j(\emptyset))$.

Then $\phi$ is equivalent to the acyclic formula
$$(\exists f)F_n(f) \wedge \hat{\phi} (f) \wedge  f(\iota(\emptyset)) = \iota^{-d_1} (x_1) \wedge \cdots \wedge f(\iota^ n(\emptyset)) = \iota ^ {-d_n} (x_n)."$$

\section {Comments}
I think what needs to be explicated is his coding function $f$ 
which is a function that sends codes of the identity of variables to codes of their stratification type.
So first one need to think of two numerical indices for each variable, one is the identity index of it and the other is the stratification type of it, the first is a natural, the second is a negative integer (or otherwise use
natural types but with reverse stratification where the type of y in each formula $y \in x$  is one step HIGHER than the type of $x$), the identity index is just a function stipulated on variables that discriminate them from each other, i.e. all occurrences of the same variable receive the same natural index, while every two distinct variables receive different natural indices. Now the function $f$ sends iterated singletons of the empty set to iterated singletons of objects substituting the relevant variables. In other words we have the following:\bigskip

$\langle a,b \rangle \in f$ if and only if: \\
there occurs a variable symbol `$v$' in the stratified formula such that :
$a$ is the $i^{th}$ iterated singleton of $ \emptyset$,  and   $b$  is  the  $d^{th}$ iterated singleton of $v$,

where $i$ is the identity index of `$v$' and $-d$ is the stratification type of `$v$'.\bigskip

Example: the formula $x \in y \wedge z \in y$,  now we'll let the identity index send $x$ to 1, $y$ to 2, and $z$ to 3,
the stratified index sends $x$ to -2 ,$ y$ to -1 and $z$ to -2, so our $f$ function will be:
$$f= \{\langle\{\emptyset\} , \{\{x\}\}\rangle, \langle\{\{\emptyset\}\}, \{y\}\rangle , \langle\{\{\{\emptyset\}\}\}, \{\{z\}\}\rangle\}$$

Now Nathan uses this function $f$ to code the variables in the formula in terms of $f$ and iterated singletons of $\emptyset$!
Once that done, then the repetition of occurrences of iterated singletons of $\emptyset$ would be always acyclic because you can rename all bound variables in those at each occurrence.\\

Another example (asked by Thomas Forster) is about how to interpret stratified but not acyclic formulas like: 
$$x \in y \wedge z \in y \wedge w \in x \wedge w \in z$$ 
Now I'll take the discriminative (identity) assignments of $x,y,z,w$ to be
1,2,3,4 respectively; while the stratification types of $x,y,z,w$ to be
-2,-1,-2,-3 respectively. So the $f$ function would be:
$$f= \{\langle\{\emptyset\},\{\{x\}\}\rangle , \langle\{\{\emptyset\}\}, \{y\}\rangle , \langle\{\{\{\emptyset\}\}\}, \{\{z\}\}\rangle, \langle\{\{\{\{\emptyset\}\}\}\}, \{\{\{w\}\}\}\rangle \}$$

Now all atomic formulas of this formula are membership atomic formulas, and those need to be interpreted in terms of $f$ images of iterated singletons of $\emptyset$. But before we go to that I'll give first the bottom line of the argument for membership, but I'll go away from the above example to a simpler one:\\

take any variables $r, s$, now take the $i+1^{th}$ iterated singleton of $r$ to be the singleton of singleton of... ($i+1$ many times iterated) of $r$, and also the $i^{th}$ iterated singleton of $s$ to be the singleton of singleton of... ($i$ many times iterated) of s.\\ 

Now examine the following formula:\\

$\exists m \forall k [k$ is the $i+1^{th}$ iterated singleton of  r $\lor k$  is  the  $i^{th}$  iterated  singleton  of  s  $\to (m$  is  an $ i+1^{th}$  iterated  element  of $k$)]$...........FORMULA(1)$\bigskip

where $i+1^{th}$ iterated element means ``element of element of element of... (iterated $i+1$ many times) of $k$"\bigskip

Now we'll see that the above formula corresponds to ``$r \in s$", it is EQUIVALENT TO ``$r \in s$"!!! [given the weak assumptions of Nathan].\bigskip

So if we denote formula(1) by $r \in_i s$, then we do have: $r\in_i s \iff r \in s$\bigskip

So this is the bottom line of the membership atomic formulas, however, to achieve acyclicity we need to code that in terms of iterated singletons of $\emptyset$ instead of using the variable symbols ``$r,s"$ in the formula! And here where the $f$ function will trip in!\bigskip

Now suppose that the discriminative index of $r$ is 1 and of $s$ is 2, and the stratification types of them are -2 of $r$ and -1 of $s$, then our function $f$ would be
$$f=\{\langle\{\emptyset\}, \{\{r\}\}\rangle, \langle\{\{\emptyset\}\}, \{s\}\rangle\}$$

Now we write up FORMULA(1) in terms of $f$, we'll have:
$$\exists m \forall k [k=f(\{\emptyset\}) \ \lor \ k=f(\{\{\emptyset\}\}) \to  m \ is \ a \ 2^{nd} \ iterated \  element \ of \ k]...FORMULA(2)$$

(where $2^{nd}$ iterated element of $k$ means an element of an element of $k$)\bigskip

Notice that this is EQUIVALENT to $``r \in s"$!\bigskip

and this would be written as: $f(\{0\}) \in_2 f(\{\{0\}\})$, where ``2" is the negative of the stratification index for r\bigskip

Generally speaking for any arbitrary variables $r$,$s$ having discriminative indices $i$,$j$ respectively, and if the stratification type of $r$ is $-d$, then we can use the following formula(2):

$$\exists m \forall k [k=f(\{..\{\emptyset\}.._i\}) \ \lor \ k=f(\{..\{\emptyset\}.._j\}) \to  m \ is \ a \ d^{th} \ iterated \  element \ of \ k]...FORMULA(2)$$

and this can be denoted as: $$f(\iota^i(\emptyset)) \in_d f(\iota^j(\emptyset))$$, which is an alternative acyclic expression of $r \in s$ \bigskip

However notice also that you only have ONE free variable in FORMULA(2), namely ``$f$", all the other symbols are either constants or are bound, and ALL of those can be renamed altogether. So suppose you have the formula $r \in s$ occurring twice which is CYCLIC but yet STRATIFIED (as you demanded), then still this can be interpreted along the above method in such a manner that iteration won't cause cyclicity because we'll rename ALL bound variables (and constants) in the above formula(2) versions of them, so the only recurrent variable would be $f$ and this will not be connected to the same variables in both occurrences (because all other variables in these formulas are re-named at each occurrence!). This is the general tactic that Nathan adopted to break cyclicity.\bigskip

Now lets discuss your formula: $``x \in y \wedge z \in y \wedge w \in x \wedge w \in z"$, the cause of the acyclicity here is clearly $w$. Now we, of course, replace each of the four atomic formulas with the formula(2) version of it, but notice that at each occurrence we'll rename all bound variables! and all constants also, and this will evaporate any cyclicity that could occur, the point is that every variable (i.e. $x$ or $y$ or $z$ or $w$) even though it recurs in the formula yet it is represented by different variables, only symbol ``$f$" is recurrent in its representation but still that won't cause any harm because, at each occurrence of ``$f$", ``$f$" is connected to a different set of variables, so no cyclicity will be raised! \bigskip

Similarly and more obviously we can write the atomic formula $r=s$
in this method, this would be:\bigskip

 $f$((discriminative identity code of r)$^{th}$ iterated singleton of $\emptyset$)=$f$((discriminative identity code of s)$^{th}$ iterated singleton of $\emptyset$)\bigskip
 
 Or using Bowler's notation, this is: $$f(\iota^i(\emptyset)) = f(\iota^j(\emptyset))$$, where $i,j$ are the identity indices of $r$ and $s$ respectively.\bigskip

Since $r$ and $s$ will receive the same stratification type, then clearly the $f$ images  will have the same degree of iteration of the ``singleton of" operator on $r$ and $s$, and clearly if the $i^{th}$ iterated singletons of $s$ and $r$ are equal then this is equivalent to saying that $s$ and $r$ are equal.\bigskip

One minor issue: I think we need to add the stipulation that if the stratified formula $\phi$ has a total of $n$ (of the metalanguage) many distinct variables, then the size of $f$ must not be more than $n$, there is an acyclic formula that expresses that:
$$\exists m_1,...,m_n \forall p \in f (p=m_1 \lor ...\lor p=m_n)$$
Now I think this can replace $F_n(f)$ big formula of Nathan, (in the final expression of the equalizing acyclic formula (the last of formulas in Nathan's proof)), so it both shortens the expression and also prevents having two distinct elements in the domain of $f$ that are $i^{th}$ singleton iterations of empty objects (for the same $i$), The problem is that when we don't have Extensionality we can have many of those in the domain of $f$, and it is better to prevent that!\bigskip

I think that Nathan's formulation will break down if there is no Extensionality over empty objects, the way how he wrote it makes it possible to have distinct same-degree iterated singletons of empty objects, and this can be disastrous.
For example, let's take the atomic formula $x=y$, let the discriminative indexing be 1 for $x$ and 2 for $y$, now let the stratification types be -1 for each of $x,y$. Now the way how Nathan wrote his formulation makes the following function $f$ possible

$$f=\{\langle \{\emptyset \}, \{x\}\rangle,  \langle \{\{\emptyset \}\}, \{y\}\rangle , \langle \{\emptyset^*\}, \{z\}\rangle \}$$ where $\emptyset$ and $\emptyset^*$ are distinct empty objects.

Now even if we have $x=y$, now if $z\neq y$, then we'll conclude that $f(\{\emptyset \}) \neq f(\{\{\emptyset \}\})$,
the reason is because the formulation doesn't discriminate between $\emptyset$ and $\emptyset^*$, so there will not exist an object $k$ that is the second projection of all elements of $f$ that has their first projections being either a set of an empty object or a set of a set of an empty object. So his method will break.

The above argument applies when we consider all $\emptyset$ symbols appearing in Nathan's $F_n(f)$ formula as parameters. But if we read  $\iota^n (\emptyset)$ as an $n$-iterated singleton of an "empty" object, then Nathan's proof is saved, albeit the formula is too long.

This will be resolved if we add the size criterion over $f$, as I've mentioned above, this size criterion can simply replace $ F_n(f)$ in his last formula, actually, the size criterion formula is much shorter than the formula $ F_n(f)$ of Nathan's.\bigskip

Requirements for this method to work are precisely axioms of : Empty, Singletons, and Boolean union. \bigskip

The use of prenex normal form turns out to be not essential for this development [Bowler], this is his take on this:\bigskip

``The trick of first reducing to prenex normal form isn't really necessary, I think, but it greatly simplifies the argument. If we want to avoid it then we have to introduce a new variable for a coding function for every quantifier and so the argument gets a bit more convoluted. Essentially if there are already variables $x_1 ... x_n$ in play, as coded by a function $f$, and we want to add a clause with an existentially quantified further variable $x_{n+1}$ then we would do so with the formula ($\exists f')F_{n+1}(f')$ and $f \subseteq f'$ and [stuff for the rest of the clause]. Universal quantifiers can be dealt with similarly. But writing the argument formally would be messy."\bigskip

To further clarify this, I'd say if we have a formula $\phi$ then we'd number its sub-formulas that has prefix quantifier strings beginning from the last and ending in $\phi$ as $\psi_1,..,\psi_n$ where $\psi_n = \phi$, now after the prefix of each $\psi_i$ we'll introduce the existence of a new coding function $f^i$ (i.e. write $(\exists f^i)$ right after the prefix of $\psi_i$) and after each $(\exists f^i)$  we'll write the last formula of Nathan's proof in terms of $f^{k \leq i}$ in the following manner:\bigskip

0. The discriminative indexing and stratification typing is stipulated over the whole formula $\phi$.\bigskip

1. The last string in Nathan's last formula, which we'll denote as the ``identification string" (the string that specifies the coding of the iterated singleton of each variable in terms of an iterated singleton of an empty set), of each function $f^i$ would only be wrote for the variables present in the prefix quantification string of $\psi_i$. When $i=n$ the identification string of $f^n$ must in addition include the coding string of all free variables in $\phi$\bigskip

2. Take the immediately preceding sub-formula $\psi_j$ that has $\psi_i$ as a sub-formula of, and add all elements of $f^j$ to elements of $f^i$ by stipulating the formula $f^j \subset f^i$ after $(\exists f^i)$.  \bigskip

3. The ``size string" of the interpreting formula (which should replace $F_n(f)$) must be equal to the sum of the size string of the preceding $\psi_j$ formula (or to the total number of free variables in $\phi$ if $i=n$) and the total number of variables in the prefix of $\psi_i$.\bigskip

4. The acyclic interpreting formula $\hat{\psi_j}$ is written in terms of $f^j$ only for the sub-formulas of $\psi_j$ outside of $\psi_i$, now $\hat{\psi_i}$ will be a sub-formula of $\hat{\psi_j}$ so all of what is in it is written in terms of $f^{k\leq i}$.\bigskip

There is another way, that is to add (to step 1) a $k^{th}$ iterated singleton of $\emptyset$ for each variable symbol  $x_k$ occurring in $\psi_i$ that is not present among variables in the prefix of $\psi_i$, to the domain of $f^i$ (by adding the formula $\exists p \in f^i (\iota^k(\emptyset)=\pi_1(p))$) and then specifying a clause for each $f^i$ to the effect that what it shares in its domain with the domain of $f^j$ (where $f^j$ immediately precede $f^i$) is to receive the same value by both functions, formally this is:
$$(\forall x) (\exists y) (\forall p) [(p \in f^i \lor p \in f^j) \wedge x=\pi_1(p)] \to y=\pi_2(p)$$. Then we proceed with the same step in 4, but step 3 need not be adopted, the size formula of each function $f^i$ will simply be equal to the total number of variables in $\psi_i$\bigskip

One last issue to be noted is that the proof is largely impredicative, since the coding function symbol `$f$' in the acyclic interpreting formula is at least one type higher than the symbol `$y$' in $\phi(y)$, so this method doesn't prove that every predicative stratified formula is equivalent to some predicative acyclic formula [Although it is easy to prove that Predicative NF is equivalent to Predicative Acyclic Comprehension]. Not only that, it can't even prove that every one degree impredicative stratified formula (where a variable symbol in $\phi(y)$ is allowed to be up to ONE level higher than $y$) is equivalent to some one-degree impredicative acyclic formula, and so whether Impedicative NF is equivalent to Impredicative Acyclic Comprehension remains unresolved.\bigskip

\noindent
Late Notes: May 12, 2020
\\
\noindent
There is something nice that I found really. I realized that we can define membership and equality actually after exactly the same acyclic formula! Now Nathan translates $x_i = x_j$  as: $$f(\iota^i (\emptyset)) = f(\iota^j(\emptyset))$$; and 
$x_i \in x_j$ as: $$f(\iota^i(\emptyset)) \in_{d} f(\iota^j(\emptyset))$$. In reality there is no need for such separate treatments, because both expressions can be defined by exactly the same acyclic formula!!!!!\bigskip

\noindent
so we translate both expressions to:  $$f(\iota^i (\emptyset)) \ \mathcal{R}_d \ f(\iota^j(\emptyset))$$; and I prefer to denote it by an undirected edge $$f(\iota^i (\emptyset)) \ -_d \ f(\iota^j(\emptyset))$$. The formula is:$$\exists e \forall p \forall k ([p \in f \land (\pi_1(p) = \iota^i(\emptyset) \lor  \pi_1(p) = \iota^j(\emptyset)) \land k=\pi_2 (p)] \to  e \in^d k)$$; where $\in^d$ is the ``$d^{th}$ iterated element of" operator.

\section{References}

[1] Al-Johar, Z; Holmes M.R., Acyclic Comprehension is equal to Stratified Comprehension, Preprint 2011. http://zaljohar.tripod.com/acycliccomp.pdf

\noindent
[2] Al-Johar, Zuhair; Holmes, M. Randall; and Bowler, Nathan. (2014). ``The Axiom Scheme of Acyclic Comprehension". Notre Dame Journal of Formal Logic, 55(1), 11-24. http://dx.doi.org/10.1215/00294527-2377851

\noindent
[3] Quine, William v. O., ``New Foundations for Mathematical Logic,” 

American Mathematical Monthly, 44 (1937), pp. 70–80.

\end{document}